\documentclass[lettersize,journal]{IEEEtran}
\usepackage{color,array,amsthm}
\usepackage{graphicx}
\usepackage{amsmath}
\usepackage{amssymb}
\usepackage{amsfonts}
\usepackage{bbm}      
\usepackage{dsfont}   
\usepackage{graphicx} 
\usepackage{amsmath,amssymb}
\usepackage{algorithm}
\usepackage{algorithmic}
\usepackage{subfig}    
\usepackage{graphicx}  
\usepackage{float}     
\newtheorem{theorem}{Theorem}

\newtheorem{problem}{Problem}
\newtheorem{remark}{Remark}
\newtheorem{assumption}{Assumption}

\newtheorem{definition}{Definition}

\begin{document}

	\title{Relaxed Control with Entropy Regularization for Itô Stochastic Systems with Input Delay} 
	
	\author{ Ruixue Li$^{1}$, Xun Li$^{2}$, Zhaorong Zhang$^{3}$, and Juanjuan Xu$^{1}$	
		\thanks{*This work was supported by the National Natural Science Foundation of China under Grants 62573262, 62503289, and the Natural Science Foundation of Shandong Province under Grant ZR2021JQ24. (Corresponding author: Juanjuan Xu.)}
		\thanks{$^{1}$Ruixue Li and Juanjuan Xu are with the School of Control Science and Engineering, Shandong University, Jinan 250061, China  {\tt\small ruixueli777@163.com; juanjuanxu@sdu.edu.cn. }}
		\thanks{$^{2}$Xun Li is with the Department of Applied Mathematics, The Hong Kong Polytechnic University, Hong Kong SAR, China {\tt\small li.xun@polyu.edu.hk.} }
		\thanks{$^{3}$Zhaorong Zhang is with the School of Computer Science and Technology, Shandong University, Qingdao 266237, China {\tt\small zhangzr@sdu.edu.cn.}}
	}
	
	\maketitle

	
	\begin{abstract}	
		This paper investigates the infinite-horizon classical stochastic optimal control problem with input delay under an entropy-regularized relaxed control framework.
		In particular, by constructing a relaxed system and introducing an entropy regularization term, we reformulate the classical optimal control problem into an entropy-regularized formulation, and derive the optimal controller that follows a Gaussian distribution.
		Furthermore, we show that the optimal Gaussian control distribution converges to the optimal Dirac measure as the exploration weight tends to zero.
		Numerical simulation is provided to validate the effectiveness of the proposed method.
	\end{abstract}
	

	\section{Introduction}\label{sec1}

	Delay is ubiquitous in control systems including finance \cite{finance},  communication networks \cite{network} and power systems \cite{power system}. It is well known that input delay can destabilize a system and degrade control performance \cite{stable1,stalbe2}, making delay compensation a central topic in control theory for decades. 
	\cite{delay1} presented a reduction method to transform linear systems with delayed control into equivalent delay-free systems for performance analysis.
	\cite{delay4} proposed a smooth adaptive state feedback controller based on a Lyapunov function to solve the tracking control problem for a class of nonlinear systems with unknown nonsymmetric dead-zone input and time delays.
	\cite{delay2} addressed prescribed-time stabilization of linear time-invariant systems with constant input delay by modeling the delay as a transport partial differential equation (PDE).
	\cite{delay3} presented a model reduction method for discrete-time coupled systems with input delays using bivariate fundamental matrices.

	The difficulty is further amplified when the system with input delay is subject to stochastic disturbances, especially multiplicative noise in the Itô stochastic systems \cite{finance,xie,wangwenjing}. 
	\cite{stoch1} established the dynamic programming principle for optimal control problems governed by stochastic differential equations with delay, derived the corresponding Hamilton–Jacobi–Bellman (HJB) equation, and analyzed the weak uniqueness of such delay systems. 
	\cite{stoch2} investigated mean-field-type stochastic optimal control problems with delay by establishing two sufficient and one necessary maximum principles for mean-field jump–diffusion stochastic delay differential equations, and applied the results to solve a delayed bicriteria mean–variance portfolio selection problem with explicit solutions.
	\cite{stoch3} studied linear quadratic regulation and stabilization problems for Itô stochastic systems with input delay via modified Riccati equations.
	\cite{stoch4} addressed the stabilization of first-order hyperbolic PDE systems with input delay and multiplicative noise by developing predictor-feedback controllers.
	All the aforementioned achievements are established within the framework of classical stochastic optimal control theory.

	In recent years, with the rapid advancement of intelligent optimization and data-driven control, reinforcement learning (RL) has attracted extensive attention.
	Among various RL framework, entropy-regularized methods explicitly incorporate exploration as a regularization term into the optimization objective and introduce a trade-off weight to adjust the policy entropy, thereby achieving a principled balance between exploration and exploitation \cite{Wang2019, Wang2020}.
	\cite{rl1} addressed the collapse of the conventional Q-function in continuous-time entropy-regularized exploratory diffusion Q-learning by introducing a first-order approximated q-function and characterizing the associated value function through martingale conditions.
	\cite{rl2} studied the exploration-exploitation trade-off and policy bias problem in continuous-time linear-quadratic RL using entropy-regularized relaxed stochastic control.
	\cite{rl3} investigated the optimal control problem for entropy-regularized backward stochastic systems by establishing a stochastic maximum principle via convex variation.
	%
	However, these existing results fail to hold for Itô stochastic systems with input delay.

	In this paper, we present the entropy-regularized relaxed control framework for Itô stochastic systems with input delay.
	The main contributions are summarized as follows.
	We first formulate the entropy-regularized stochastic optimal control problem for delayed systems, and prove that the resulting optimal policy obeys a Gaussian distribution. 
	Moreover, we prove that the optimal Gaussian control distribution converges to the optimal Dirac measure as the exploration weight tends to zero.
	Numerical simulation is carried out to verify the effectiveness of the proposed approach.
	
	The remainder of this paper is organized as follows.
	Section \ref{sec2} formulates the infinite-horizon classical stochastic optimal control problem with input delay.
	Section \ref{sec3} develops the entropy-regularized relaxed control framework and derives the optimal Gaussian control distribution.
	Section \ref{sec5} provides a numerical simulation to demonstrate the effectiveness of the proposed method.

	The notations adopted throughout this paper are defined as follows. $\mathbb{R}^n$ stands for the set of all $n$-dimensional real vectors, and $\mathbb{R}^{n \times m}$ represents the set of $n\times m$-dimensional real matrices. 
	Matrix $M>0$ $(\ge 0)$ is strictly positive-definite (positive semi-definite).
	$H'$ denotes the transpose of  matrix $H$.
	$\bar{C}[ -h,0 )$ denotes the space of continuous bounded functions mapping $[ -h,0 ) $ to $ \mathbb{R}^m $, which is formally defined as $\bar{C}[ -h,0 )=\{\alpha(t):[-h,0) \to \mathbb{R}^m \ is \ continuous \ and \ \sup_{-h\le t< 0} \left \| \alpha(t) \right \| < \infty  \}$.	
	Let $\{\Theta ,\mathcal{F},\mathbb{P},\{\mathcal{F}_s\}_{s\ge0}\}$ be a complete stochastic basis, where $\mathcal{F}_0$ contains all P-null sets in $\mathcal{F}$, and the filtration $\{\mathcal{F}_s\}_{s\ge0}$ is generated by the standard Brownian motion $\{\omega(s)\}_{s\ge0}$.
	The conditional expectation of $x(t)$ with respect to $\mathcal{F}_{s}$ is defined as $\hat{x}(t\mid s) \doteq E[x(t) \mid \mathcal{F}_{s}]$ .

	\section{Problem formulation}\label{sec2}
	In this section, we formulate the infinite-horizon classical stochastic optimal control problem with input delay. The considered system is governed by the following Itô stochastic differential equation involving multiplicative noise and delayed control input:
	\begin{align}	\label{classical system}
		\begin{cases} 
			dx(t)\hspace{-2pt}=\hspace{-2pt} [Ax(t)\hspace{-2pt}+\hspace{-2pt}Bu(t-h)]dt\hspace{-2pt}+\hspace{-2pt} [\bar{A}x(t)\hspace{-2pt}+\hspace{-2pt}\bar{B}u(t-h)]d\omega(t),   \\
			x(0)=x, u(\tau)=\alpha(\tau), \tau \in [ -h,0 ),
		\end{cases}
	\end{align}	
	where $x(t)\in \mathbb{R}$ is the state, $u(t)\in \mathbb{R}^m $ is the control input, and $\omega(t)$ is the standard Brownian motion. 
	The constant $h>0$ denotes the input delay. $x$ and $\alpha(\tau) \in \bar{C}[ -h,0 )$ are initial values and deterministic for notational simplicity. $A, \bar{A}, B, \bar{B}$ are constant matrices with appropriately compatible dimensions.
	
	The discounted infinite-horizon cost function is
	\begin{align}\label{classical cost}
		J=& \frac{1}{2}  E \{  \int_{0}^{\infty }   e^{-\rho t} [x'(t)Qx(t) +  u'(t-h)Ru(t-h)]dt   \} ,
	\end{align}	
	where  $\rho > 0$ is the discount rate, $Q$ and $R$ are symmetric matrices with compatible dimensions.
	
	The following standard assumption is made for analyzing the system stabilization.
	
	\begin{assumption}\label{assump1}
		Matrices $Q=C'C$ and $R$ are both positive definite.
	\end{assumption}
	
	
	\begin{definition}
		Denote by $\mathcal{U}(x)$ the set of all admissible controls,
		$U\doteq  \{u(t), t \in [-h,\infty) \}\in \mathcal{U}(x)$ if 
		\begin{enumerate}
			\item $u(t)$ is $\mathcal{F}_{t}$-adapted;
			\item For each $t \ge 0$, $E\{\int_0^t u'(s)u(s) ds\} < \infty$;
			\item With $x(t)$ solving \eqref{classical system},
			$\lim_{t\to\infty} e^{-\rho t} E\left[x'(t)x(t)\right] = 0;
			$
			\item $  J < \infty.$
		\end{enumerate}
	\end{definition}

	\begin{problem}\label{problem 1}
		Find an admissible control
		$U \in \mathcal{U}(x) $ with constant matrix gain that minimizes the cost function \eqref{classical cost} subject to the system \eqref{classical system}.	
	\end{problem}

	\begin{remark}
		When $\rho=0$, the solution to Problem \ref{problem 1} has been derived in \cite{stoch3}, which is within the classical stochastic control framework. As entropy regularization and relaxed control theory advances, it is necessary to further investigate this problem via an entropy-regularized relaxed control framework, which will be presented in the next section.
	\end{remark}

	\section{Entropy-regularized stochastic control problem}\label{sec3}
	In this section, we introduce an entropy-regularized relaxed control framework for continuous-time Itô stochastic systems with input delay. Motivated by \cite{Wang2019}, we reformulate the original system dynamics \eqref{classical system} into an exploratory formulation, 
	which is modelled by a distribution $\pi \doteq \{ \pi_t, t\in [-h,\infty) \}$.
	The resulting relaxed system is described as:	
	\begin{align}	\label{E system}
		\begin{cases} 
			dX^{\pi}_{t} = [AX^{\pi}_{t}+B\mu_{t-h}]dt +  D (X^{\pi}_{t},{\pi}_{t-h}) d\omega(t),   \\
			X^{\pi}_{0}=x, \pi_{\tau}(u) = \delta_{\alpha(t)} , \tau\in[ -h,0 ) ,
		\end{cases}	
	\end{align}
	where $\delta_{\alpha(t)}$ is a Dirac distribution, and the diffusion term $ D (X^{\pi}_{t},{\pi}_{t-h})$ satisfies
	\begin{align}
		D (X^{\pi}_{t},{\pi}_{t-h}) D' &(X^{\pi}_{t},{\pi}_{t-h})=  \int_{\mathbb{R}^m}    [\bar{A}X^{\pi}_{t}+\bar{B}u]\nonumber\\
		& \times[\bar{A}X^{\pi}_{t}+\bar{B}u]' \pi_{t-h}(u)du,
	\end{align}
	while the mean and variance are defined as
	\begin{align}
		&	\mu_{t-h}= \int_{\mathbb{R}^m} u\pi_{t-h}(u)du,\\
		&	\Sigma_{t-h} = \int_{\mathbb{R}^m} uu'\pi_{t-h}(u)du-\mu_{t-h}\mu_{t-h}'.
	\end{align}

	To quantify the overall exploration level, we adopt the Shanon's differential entropy 
	\begin{align}
		-\int_{\mathbb{R}^m} \pi_{t-h}(u) \ln \pi_{t-h}(u) du  
	\end{align}
	as the exploration metric and introduce the exploration weight $\lambda \ge 0$ to balance exploration and exploitation. Combining with the classical cost function of infinite-horizon cost function \eqref{classical cost}, we construct the entropy-regularized exploratory cost function as:
	\begin{align}\label{entropy cost}
		J^{\pi}  =&   E \{    \int_{0}^{\infty }   e^{-\rho t}  \int_{\mathbb{R}^m}  [\frac{1}{2} (X^{\pi}_{t})'QX^{\pi}_{t} + \frac{1}{2}  u'Ru \nonumber \\
		& + \lambda \ln \pi_{t-h}(u)] \pi_{t-h}(u)du   dt   \}.
	\end{align}

	We now define the admissible control set $\Pi(x)$ and formulate the entropy-regularized stochastic optimal control problem with input delay.
	
	\begin{definition}
		Let $\mathcal{P}(\mathbb{R}^m)$ and $\mathcal{B}(\mathbb{R}^m)$  represent the set of all density functions and Borel algebra on $\mathbb{R}^m$. A relaxed control policy $\pi \in \Pi(x)$ if it satisfies the following conditions:
		\begin{enumerate}
			\item For each $t $, $\pi_t \in \mathcal{P}(\mathbb{R}^m)$ almost surely;
			\item For each $G \in \mathcal{B}(\mathbb{R}^m)$, $\left\{\int_G \pi_t(u)du\right\}$ is $\mathcal{F}_t$-progressively measurable;
			\item For each $t $, $E\{\int_0^t \left(\mu_s\mu_s'+ tr(\Sigma_s) \right) ds\} < \infty$;
			\item With $X^{\pi}_t$$ solving \eqref{E system}, 
			\lim_{t\to\infty} e^{-\rho t} E\left[(X_t^{\pi})'X_t^{\pi}\right] = 0;
			$
			\item $J^{\pi}< \infty$.
		\end{enumerate}
	\end{definition}

	\begin{problem}\label{problem 2}
		Find an optimal control $\pi \in \Pi(x)$ that minimizes the entropy-regularized cost function \eqref{entropy cost} subject to system \eqref{E system}.	
	\end{problem}

	
	We next present 
	the value function and the controller of the entropy-regularized stochastic optimal control problem (Problem \ref{problem 2}) and those of the classical problem (Problem \ref{problem 1}).
	To this end, we first define the value function of Problem \ref{problem 2} as
	\begin{align} \label{Vpi}
		V^{\pi}(x)=&  \min_{\pi \in \mathcal{P}( \mathbb{R}^m ) } J^{\pi},
	\end{align}
	and that of Problem \ref{problem 1} as
	\begin{align}\label{V}
		V(x)=& \min_{u }  J   .
	\end{align}
	Then, we introduce the modified algebraic Riccati equation (ARE):
	\begin{align}\label{classical RE}
		0 = Q + PA +A'P+\bar{A}'P\bar{A}-\Lambda(h)- \rho P ,
	\end{align}
	where
	\begin{align}
		&\Lambda(s) =   e^{(A-\frac{1}{2}\rho)'s}\Lambda(0)e^{(A-\frac{1}{2}\rho)s},\label{Lambdas}\\
		&\Lambda(0) =  K'\Omega K, \label{Lambda0} \\
		&K = -\Omega^{-1}[B'P+ \bar{B}'P\bar{A}-B'\int_{0}^{h}\Lambda(s)ds], \label{K} \\
		&\Omega  =  R + \bar{B}'P\bar{B}>0. \label{Omega}
	\end{align}

	\begin{theorem}\label{theorem: main}
		Under Assumption \ref{assump1}, suppose that the ARE \eqref{classical RE}--\eqref{Omega} admits a unique positive definite solution \(P-\int_{0}^{h}\Lambda(s)ds>0\), the following statements hold.
		
		\noindent\textbf{(a)} 
		\begin{align}\label{Vpifunction}		
			V^{\pi}(x)=&\frac{1}{2} \{ x'Px-x'\int_{0}^{h}\Lambda(s)\hat{x}_{0\mid s-h}ds\} \nonumber \\
			& -  \frac{\lambda }{2\rho} [ \ln((2\pi e)^m \left |  \lambda\Omega^{-1}   \right | )  -1    ],
		\end{align} is the value function of Problem \ref{problem 2}, and the corresponding optimal
		Gaussian control distribution is given by 
		\begin{align} \label{pi*}
			\pi^*_{t-h}= \mathcal{N}\left( u \;\middle|\; \mu^*_{t-h},\; \Sigma^*_{t-h} \right), 
		\end{align}
		where
		\begin{align}\label{mean}
			&	\mu^*_{t-h} = K\hat{X}^{\pi}_{t\mid t-h} ,\\
			&	\Sigma^*_{t-h} =  \lambda(R+\bar{B}' P\bar{B} )^{-1} \label{var},\\
			&	\hat{X}^{\pi}_{t\mid t-h} =e^{Ah}X^{\pi}_{t-h} +  \int_{t-h}^{t} e^{A(t-s)}B\mu_{s-h}ds . \label{Xhat}
		\end{align} 

		\noindent\textbf{(b)} 		
		\begin{align} \label{Vfunction}
			V(x)=\frac{1}{2} \{  x'Px-x'\int_{0}^{h}\Lambda(s)\hat{x}_{0\mid s-h}ds\} .
		\end{align}
		is the value function of Problem \ref{problem 1}, 
		and the corresponding optimal control is
		\begin{align}\label{u*}
			u^*(t-h) =   K \hat{x}_{t\mid t-h},
		\end{align}
		where
		\begin{align}
			\hat{x}_{t\mid t-h} = &e^{Ah}x(t-h) +  \int_{t-h}^{t} e^{A(t-s)}Bu(s-h)ds \label{xhat}.
		\end{align}	

		\noindent\textit{Proof:} The detailed proof is given in Appendix A. \hfill $\blacksquare$
	\end{theorem}

	From Theorem \ref{theorem: main}, the optimal control of Problem \ref{problem 2} obeys a Gaussian distribution with mean \(\mu^*_{t-h}=K\hat{X}^{\pi}_{t\mid t-h}\), 
	and variance \(\Sigma^*_{t-h}=\lambda\Omega^{-1}\) related to the regularization parameter \(\lambda\). As \(\lambda\) decreases, the variance of the Gaussian policy gradually decays.
	Naturally, 	
	as the regularization effect gradually diminishes, 
	the entropy-regularized problem converges to its classical counterpart,
	which is formally characterized in the following theorem.
	\begin{theorem}	\label{theorem dirac}
		As the exploration weight $\lambda$ vanishes, 
		the optimal Gaussian control distribution \eqref{pi*} converges to the optimal Dirac measure, i.e., 
		\begin{align} \label{eq:dirac}
			\lim_{\lambda  \to 0 } \pi^*_{t-h} = \delta_{u^*(t-h)},  weakly,
		\end{align}
		where $\delta_{u^*(t-h)}$ denotes the optimal Dirac measure associated with $u^*(t-h)$ given in \eqref{u*}.
		
		$\mathit{\textbf{Proof:}}  $ The proof is referred to \cite{Wang2019} and omitted herein. 	\hfill $\blacksquare$		
	\end{theorem}

	\begin{remark}
		Theorem \ref{theorem: main} establishes an entropy-regularized relaxed control framework for infinite-horizon Itô stochastic systems with input delay. 
		The resulting optimal relaxed control \eqref{pi*} follows a Gaussian distribution.
		Moreover, the optimal Gaussian control distribution converges to the optimal Dirac measure as the regularization weight decays to zero as shown in Theorem \ref{theorem dirac}.
	\end{remark}



	\section{Simulation}\label{sec5}
	In this section, we present a numerical example to validate the theoretical results and examine the effectiveness of the developed entropy-regularized relaxed control framework for stochastic systems with input delay.

	In particular, the system parameters in \eqref{classical system} are chosen as \(A = 1.2\), \(\bar{A} = 0.6\), \(B = 2\) and \(\bar{B} = 0.2\). For the cost function \eqref{classical cost}, the state weighting matrix and control weighting matrix are set to be \(Q = 1\) and \(R = 0.1\), respectively. The initial state is assigned as \(x = 30\), the input delay is selected as \(h = 1\),  and the control input satisfies \(\alpha(t) = 0\) for \(t\in[-h,0)\). Furthermore, the discount factor is set to be \(\rho = 30\) and the total simulation time is \(T = 20\).

\begin{figure}[H]
		\begin{center}
			\includegraphics[width=0.45\textwidth, keepaspectratio]{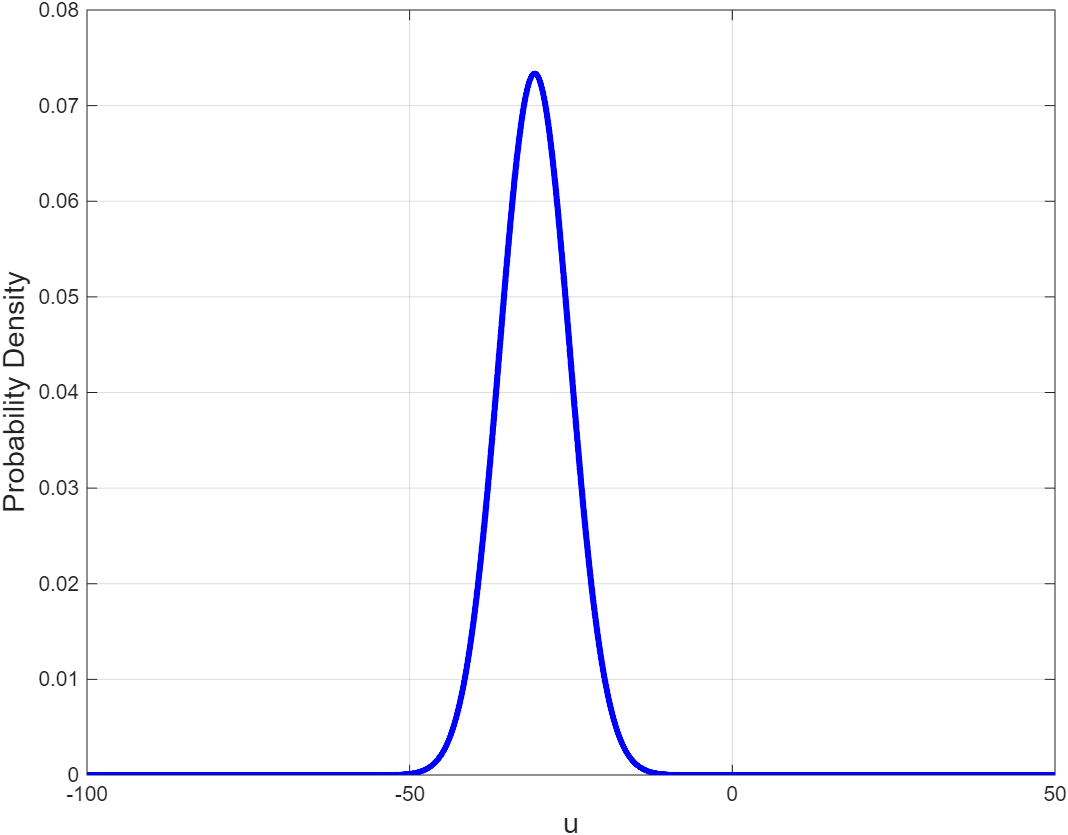}
			\caption{Gaussian distribution at t=4, $\lambda$ = 3.000.}
			\label{fig:Gaussian4}
		\end{center}
	\end{figure}

\begin{figure}[H]
		\begin{center}
			\includegraphics[width=0.45\textwidth, keepaspectratio]{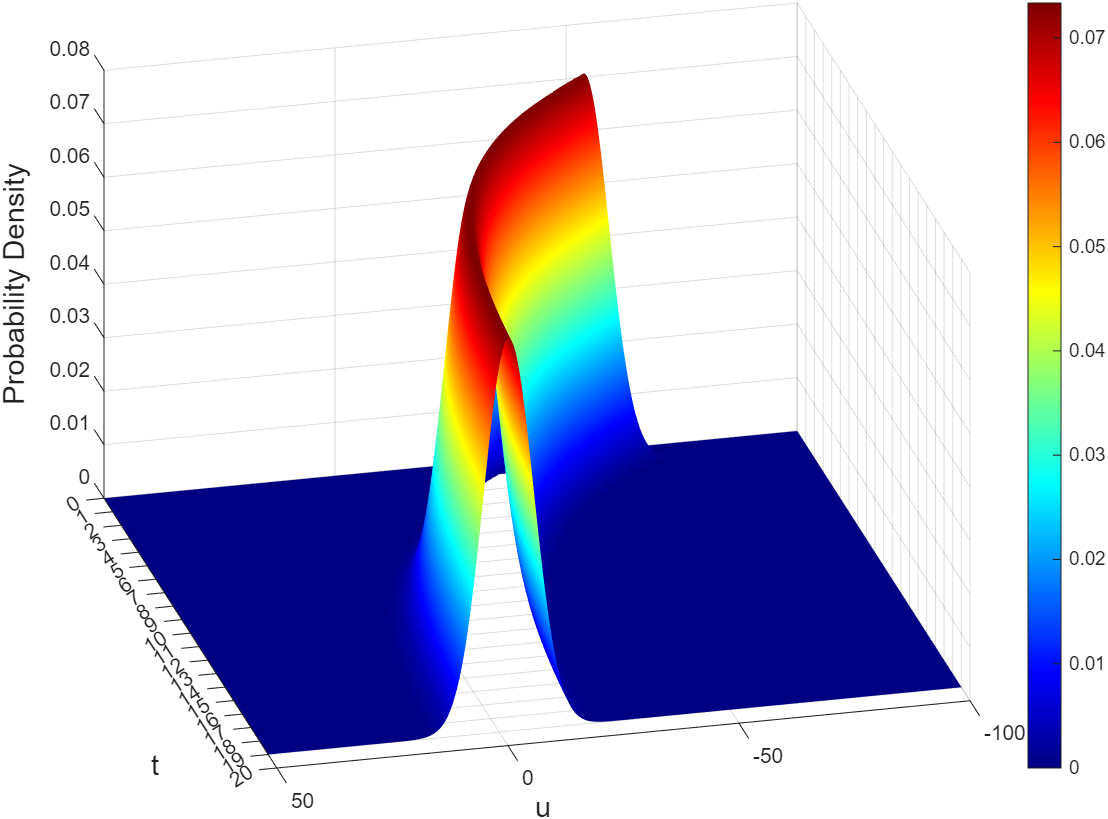}
			\caption{All-time Gaussian distribution $\lambda$ = 3.000.}
			\label{fig:Gaussian3D}
		\end{center}
	\end{figure}

\begin{figure}[H]
		\begin{center}
			\includegraphics[width=0.45\textwidth, keepaspectratio]{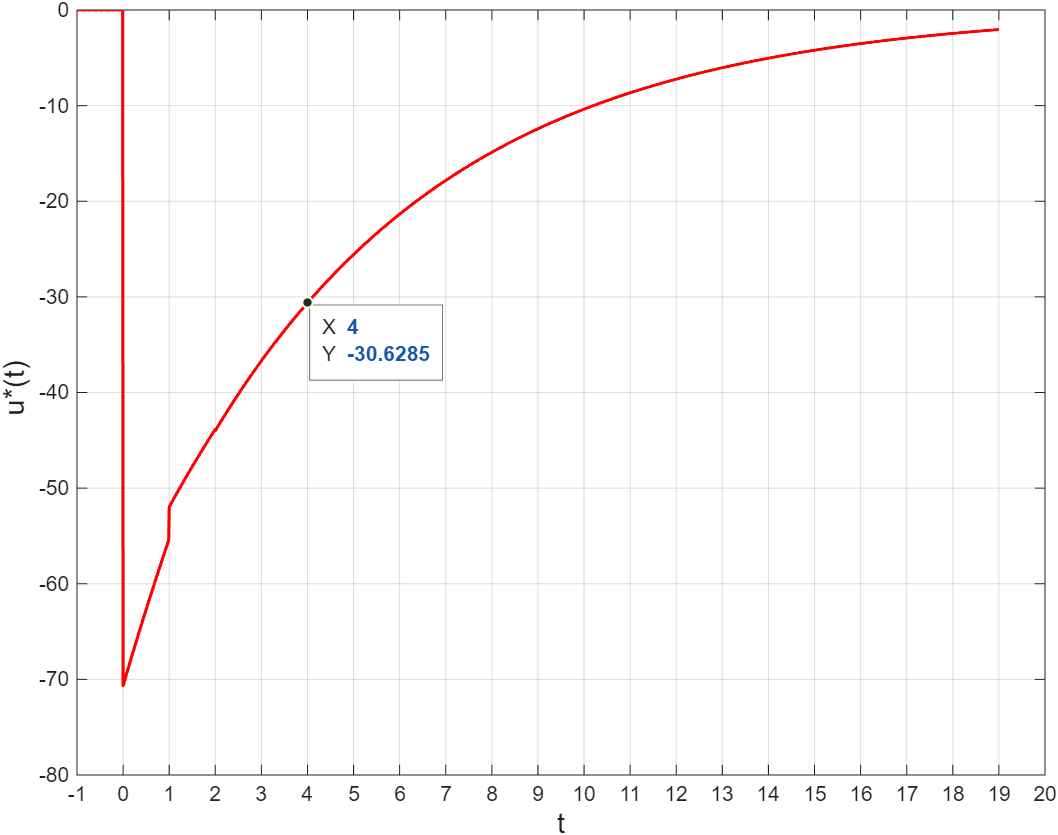}
			\caption{The trajectory of control input.}
			\label{fig:u}
		\end{center}
	\end{figure}

	We now verify the theoretical results in Theorem \ref{theorem: main}.
	By solving the modified algebraic Riccati equation \eqref{classical RE}--\eqref{Omega}, we obtain the solutions \( P = 0.036232 \) and \( K = -0.709606 \).
	Within the proposed entropy-regularized relaxed control framework, the relaxed controller derived from \eqref{pi*} follows a Gaussian distribution are shown in Fig. \ref{fig:Gaussian4} and Fig. \ref{fig:Gaussian3D}. 
	Specifically, Fig. \ref{fig:Gaussian4} shows the Gaussian distribution at \(t=4\), and Fig. \ref{fig:Gaussian3D} illustrates the time-varying Gaussian distribution over $[0,T-h]$.	
	Within the classical stochastic control framework, the classical optimal controller is further obtained via \eqref{u*}, which is depicted in Fig. \ref{fig:u}.

\begin{figure}[H]
		\begin{center}
			\includegraphics[width=0.45\textwidth, keepaspectratio]{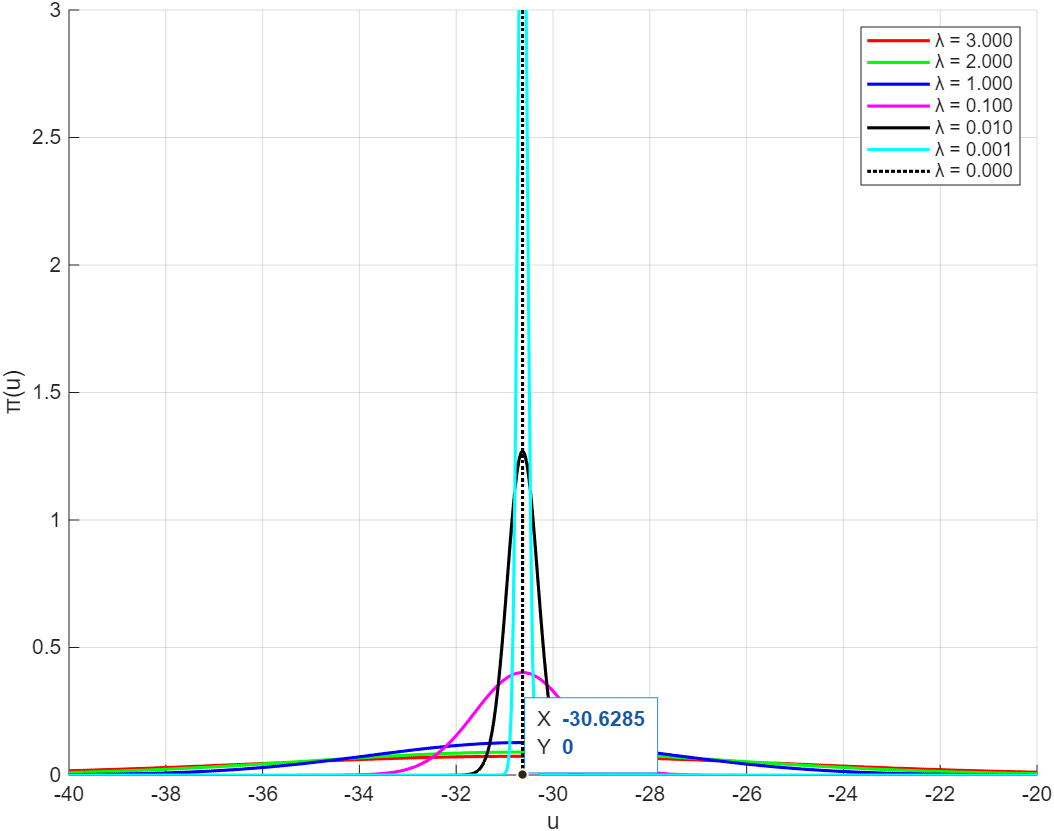}
			\caption{Gaussian distributions at t=4 under different $\lambda$.}
			\label{fig:dirac}
		\end{center}
	\end{figure}

	To further verify the convergence property of the relaxed control stated in Theorem \ref{theorem dirac}, we compute the optimal controller for different \(\lambda \in \{3, 2, 1, 0.1, 0.01, 0.001, 0\}\) via \eqref{pi*}, 
	as illustrated in Fig. \ref{fig:dirac}, which shows the resulting Gaussian distributions at \(t=4\) under different \(\lambda\).
	The marked points in Fig. \ref{fig:u}--\ref{fig:dirac} show that as \(\lambda\) decays to zero, the mean \(\mu^*_{t-h}=-30.6285\) of the optimal Gaussian distribution \(\pi^*_{t-h}\) equals to the optimal classical control \(u^*(t-h) = -30.6285\), and the distribution converges to the Dirac measure, which verifies Theorem 2.
	
	%

	\section*{APPENDIX A}
	%
	%
	%

	In this appendix, we derive the admissible optimal controllers and value functions of Problem \ref{problem 1} and \ref{problem 2} respectively.
	
	\noindent\textbf{(i)} 
	We first focus on Problem \ref{problem 1} and begin its derivation by constructing the following function:
	\begin{align} \label{y}
		y(t,x)=e^{-\rho t}  \{ x'(t)Px(t)-x'(t)\int_{t}^{t+h}\Lambda(s-t)\hat{x}_{t\mid s-h}ds\},
	\end{align}
	where \(P\), $\Lambda(\cdot)$  are shown in \eqref{classical RE} and \eqref{Lambdas}.
	
	By applying Itô’s formula to $y(t,x)$, it yields that
	\begin{align} \label{dy}
		dy(t,x) =&   e^{-\rho t}\{x'(t) [A'P+PA+\bar{A}'P\bar{A} -\rho P -  \Lambda(h)]x(t) \nonumber \\
		& +  u'(t-h)[B'P+\bar{B}'P\bar{A} ]x(t)  +  x'(t)[PB+\bar{A}'\nonumber \\
		&\times P\bar{B} ]u(t-h)  - u'(t-h)B' \int_{t}^{t+h}\Lambda(s-t) \nonumber \\
		& \times \hat{x}_{t\mid s-h}ds -x'(t) \int_{t}^{t+h}\Lambda(s-t)dsBu(t-h)   \nonumber \\  
		&+ u'(t-h)  \bar{B}'P\bar{B} u(t-h)+ x'(t)\Lambda(0)\hat{x}_{t\mid t-h}  \} dt\nonumber\\
		&  +  \star dw(t),
	\end{align}
	where \(\star dw(t)\) denotes the diffusion term, whose explicit form is irrelevant to the subsequent solution. 
	
	Note that all admissible controllers of Problem \ref{problem 1} are such that \(\lim_{t\to\infty} e^{-\rho t}E\big[x'(t)x(t)\big] = 0\).
	With initial state \(x(0)=x\), we have
	\begin{align*}
		E{\int_{0}^{\infty} dy(s,x) } =& - E\{   x'Px-x'\int_{0}^{h}\Lambda(s)\hat{x}_{0\mid s-h}ds \}.
	\end{align*}
	Taking integral from $0$ to $\infty$ on both sides of \eqref{dy} and simplifying the resulting terms via \eqref{classical RE}--\eqref{Omega} yields
	\begin{align*} \label{y=}
		-& E\{     x'Px-x'\int_{0}^{h}\Lambda(s)\hat{x}_{0\mid s-h}ds\} 
		\nonumber \\ 
		=& E\{\int_{0}^{\infty}   e^{-\rho \tau}\{x'(\tau) [A'P+PA+\bar{A}'P\bar{A} -\rho P -  \Lambda(h)]x(\tau) 
		\nonumber \\
		& +  u'(\tau-h)[B'P+\bar{B}'P\bar{A} ]x(\tau)  +  x'(\tau)[PB+\bar{A}'P\bar{B} ]
	\end{align*}
	\begin{align}
		&\times u(\tau-h) \hspace{-2pt} -\hspace{-2pt} u'(\tau-h)B' \int_{0}^{h}\Lambda(s)\hat{x}_{\tau\mid \tau+s-h}ds-x'(\tau) \nonumber \\
		&\times  \int_{0}^{h}\Lambda(s)dsBu(\tau-h)  + u'(\tau-h)  \bar{B}'P\bar{B}u(\tau-h) \nonumber \\  
		&   + x'(\tau)\Lambda(0)\hat{x}_{\tau\mid \tau-h}  \} d\tau \}
		\nonumber \\ 
		=& E\{\int_{0}^{\infty}\hspace{-3pt}   e^{-\rho \tau}\{- x'(\tau) Qx(\tau) -  u'(\tau-h)\Omega Kx(\tau) - x'(\tau)\nonumber \\  
		& \times K' \Omega u(\tau-h) \hspace{-2pt}+\hspace{-2pt} u'(\tau-h)  \bar{B}'P\bar{B}u(\tau-h)\hspace{-2pt}+\hspace{-2pt} x'(\tau)\Lambda(0)\nonumber\\
		&\times \hat{x}_{\tau\mid \tau-h}  \} d\tau \}\nonumber \\
		=& E\{\int_{0}^{\infty} \hspace{-3pt}  e^{-\rho \tau}\{- x'(\tau) Qx(\tau) -  u'(\tau-h)\Omega Kx(\tau)- x'(\tau) \nonumber \\  
		& \times K' \Omega u(\tau-h) + u'(\tau-h) \Omega u(\tau-h)+ x'(\tau)\Lambda(0) \nonumber \\  
		&\times \hat{x}_{\tau\mid \tau-h}- u'(\tau-h)R u(\tau-h)  \} d\tau \}.
	\end{align}
	In view of the definition of \(V(t)\) given in \eqref{V}, we can rewrite the value function as
	\begin{align} \label{V*}
		V(x)
		= &\min_{u} \frac{1}{2} E\{ \int_{0}^{\infty}  e^{-\rho \tau}\{ -  u'(\tau-h)\Omega Kx(\tau)+ x'(\tau)\Lambda(0) \nonumber \\  
		& \times \hat{x}_{\tau\mid \tau-h} - x'(\tau)K' \Omega u(\tau-h)+  u'(\tau-h) \Omega  \nonumber \\
		& \times u(\tau-h)\} d\tau+ [ x'Px-x'\int_{0}^{h}\Lambda(s)\hat{x}_{0\mid s-h}ds] \}
		\nonumber \\
		=&\min_{u} \frac{1}{2} E\{\int_{0}^{\infty}   e^{-\rho \tau}\{ - u'(\tau-h)\Omega K x(\tau)  - x'(\tau)K'  \nonumber\\
		&\times \Omega u(\tau-h)+ u'(\tau-h) \Omega u(\tau-h) + x'(\tau)K'\Omega K
		\nonumber \\
		&\times \hat{x}_{\tau\mid \tau-h} \} d\tau +  [ x'Px - x'\int_{0}^{h}\Lambda(s) \hat{x}_{0\mid s-h}ds] \}
		\nonumber \\
		=&\min_{u} \frac{1}{2}  E\{  \hspace{-2pt} \int_{0}^{\infty}   \hspace{-3pt} e^{-\rho \tau}[  (u(\tau-h) \hspace{-2pt} - \hspace{-2pt}K\hat{x}_{\tau\mid \tau-h})' \Omega (u(\tau-h) \nonumber \\
		&-K \hat{x}_{\tau\mid \tau-h}) ] d\tau +  [ x'Px-x' \int_{0}^{h}\Lambda(s) \hat{x}_{0\mid s-h}ds ] \} \nonumber \\
		=&  \frac{1}{2}  \{  x'Px-x' \int_{0}^{h}\Lambda(s)\hat{x}_{0\mid s-h}ds  \}.
	\end{align}
	Then, we define the Lyapunov function as \(V(x)\) in \eqref{V*}. For any \(x(t)\neq 0\), it yields that 
	\begin{align}
		V(x) &=  \frac{1}{2}  \{ x'Px-x' \int_{0}^{h}\Lambda(s)\hat{x}_{0\mid s-h}ds  \}\nonumber\\
		& \ge   \frac{1}{2}  \{  x' [P-\int_{0}^{h}\Lambda(s)ds]x \}\nonumber\\ 
		& >  0,
	\end{align}
	and
	\begin{align}
		\dot{V}(t) &= - \frac{1}{2}  \{x'(t)Qx(t)+ u'(t-h)Ru(t-h) \}\nonumber\\ 
		& <  0,
	\end{align}
	by virtue of $P-\int_{0}^{h}\Lambda(s)\mathrm{d}s>0$, $Q>0$ and $R>0$.		
	According to the Lyapunov stability theorem \cite{Lyapunov stability theorem}, it follows that $\lim_{t\to\infty} e^{-\rho t} E\left[x'(t)x(t)\right] = 0$. 
	Therefore, \eqref{u*} is the admissible optimal control of Problem \ref{problem 1}, and the corresponding value function is given by \eqref{Vfunction}.

	\noindent\textbf{(ii)} Similarly, we proceed to solve the optimal controller and value function of Problem \ref{problem 2}. 
	To this end, we introduce the function
	\begin{align}\label{ypi}		
		y^{\pi}(t,X^{\pi}_{t})=& e^{-\rho t}  \{(X^{\pi}_{t})'PX^{\pi}_{t}-(X^{\pi}_{t})'\int_{t}^{t+h}\Lambda(s-t)\nonumber \\
		&\times \hat{X}^{\pi}_{t\mid s-h}ds \},
	\end{align}
	where the matrices \(P\), $\Lambda(\cdot)$ and \(\Omega\) satisfy \eqref{classical RE}-- \eqref{Omega}.

	Applying Itô’s formula to $y^{\pi}(t,X^{\pi}_{t})$ in \eqref{ypi}, it yields that
	\begin{align} \label{dypi}
		dy^{\pi}(t,X^{\pi}_{t})
		=&  -\rho   e^{-\rho t} \{(X^{\pi}_{t})'PX^{\pi}_{t}-(X^{\pi}_{t})'\int_{t}^{t+h}\Lambda(s-t) \nonumber\\
		& \times \hat{X}^{\pi}_{t\mid s-h}ds \} + e^{-\rho t} \{  (AX^{\pi}_{t} + B \mu_{t-h} )'P X^{\pi}_{t}   \nonumber \\
		& +  (X^{\pi}_{t})'P (AX^{\pi}_{t} + B \mu_{t-h} ) +     D' (X^{\pi}_{t},{\pi}_{t-h})  P  \nonumber \\
		&\times  D (X^{\pi}_{t},{\pi}_{t-h})- (X^{\pi}_{t})'\Lambda(h) X^{\pi}_{t}  +  (X^{\pi}_{t})'\Lambda(0)\nonumber\\
		& \times \hat{X}^{\pi}_{t\mid t-h}   -    (AX^{\pi}_{t} + B \mu_{t-h} )'  \int_{t}^{t+h}\Lambda(s-t) \nonumber \\
		&\times \hat{X}^{\pi}_{t\mid s-h}ds   \hspace{-2pt}- \hspace{-2pt} (X^{\pi}_{t})'\hspace{-2pt} \int_{t}^{t+h} \hspace{-4pt} \Lambda(s-t)   ( A \hat{X}^{\pi}_{t\mid s-h}\nonumber\\
		& +B \mu_{t-h}  )  ds - (X^{\pi}_{t})' \int_{t}^{t+h}  \frac{\partial  \Lambda(s-t)    }{\partial t} \nonumber\\
		&\times  \hat{X}^{\pi}_{t\mid s-h} ds  \} dt + * dw(t)		
		\nonumber \\
		=&   e^{-\rho t} \{ \int_{\mathbb{R}^m}  \{ (X^{\pi}_{t})' [A'P+PA+\bar{A}'P\bar{A} -\rho P  \nonumber \\
		& \hspace{-2pt} -  \Lambda(h)]X^{\pi}_{t}  \hspace{-2pt} + \hspace{-2pt}  u'[B'P+\bar{B}'P\bar{A} ]X^{\pi}_{t}   \hspace{-2pt}+ \hspace{-2pt}  (X^{\pi}_{t})'[PB \nonumber \\
		&  \hspace{-2pt} +\bar{A}'P\bar{B} ]u  \hspace{-2pt} - \hspace{-2pt} u'B'\hspace{-2pt} \int_{t}^{t+h}\hspace{-2pt}\Lambda(s-t)\hat{X}^{\pi}_{t\mid s-h}ds \hspace{-2pt}   \nonumber \\
		&  -\hspace{-2pt}(X^{\pi}_{t})'\int_{t}^{t+h}\Lambda(s-t)dsBu+ u'  \bar{B}' P\bar{B}u  \nonumber \\
		& + (X^{\pi}_{t})'\Lambda(0) \hat{X}^{\pi}_{t\mid t-h}  \} \pi_{\tau-h}(u)du\} dt \nonumber\\
		& +  * dw(t),
	\end{align}
	where $*dw(t)$ denotes the diffusion term, which is omitted since it plays no role in the following derivation.
	
	Integrating both sides of \eqref{dypi} from 0 to \(\infty\) and simplifying the resulting expression by virtue of \eqref{classical RE}--\eqref{Omega}, we arrive at
	\begin{align} \label{ypii1}
		E\{ \int_{0}^{\infty} \hspace{-2pt}	 dy^{\pi}(t,X^{\pi}_{t} ) \} 
		=& E\{  \int_{0}^{\infty}  \hspace{-2pt} e^{-\rho t}  \hspace{-2pt} \int_{\mathbb{R}^m}  \{ - (X^{\pi}_{t})' QX^{\pi}_{t}-u'R \nonumber\\
		&  \times u +   (u- K \hat{X}^{\pi}_{t\mid t-h}  )'  \Omega (u- K  \nonumber \\
		& \times  \hat{X}^{\pi}_{t\mid t-h}  )  \}\pi_{\tau-h}(u)du  dt \}.
	\end{align}

	For all admissible controllers of Problem \ref{problem 2}, we have \(\lim_{t\to\infty} e^{-\rho t}E\big[(X^{\pi}_{t})'X^{\pi}_{t}\big] = 0\).
	With the initial state \(X^{\pi}_{0}=x\), we thus conclude that
	\begin{align}\label{ypii2}
		E\{\int_{0}^{\infty} dy^{\pi}(t,X^{\pi}_{t} ) \} =& - E\{   x'Px-x'\int_{0}^{h}\Lambda(s)\hat{x}_{0\mid s-h}ds \}.
	\end{align}
	
	Combining the above results \eqref{ypii1}--\eqref{ypii2} with the value function definition \eqref{Vpi}, we rearrange the value function of Problem \ref{problem 2} into the following form:
	\begin{align} \label{reypi}
		V&^{\pi}(x)\nonumber\\
		=& \min_{\pi \in \mathcal{P}( \mathbb{R}^m ) }   E \{  \int_{0}^{\infty }   e^{-\rho t} \int_{\mathbb{R}^m}  [\frac{1}{2} (X^{\pi}_{t})'QX^{\pi}_{t} + \frac{1}{2} u'R u \nonumber \\
		&  + \lambda \ln \pi_{t-h}(u)] \pi_{t-h}(u)du dt  \}\nonumber\\
		=& \min_{\pi \in \mathcal{P}( \mathbb{R}^m ) }   E \{  \int_{0}^{\infty }   e^{-\rho t}\hspace{-4pt} \int_{\mathbb{R}^m}  [\frac{1}{2}(X^{\pi}_{t})'QX^{\pi}_{t}+  \frac{1}{2} u'Ru]  \nonumber \\
		&\times  \pi_{t-h}(u) dudt\} +E \{  \int_{0}^{\infty } e^{-\rho t} \int_{\mathbb{R}^m} \lambda\ln \pi_{t-h}(u)  \nonumber \\
		& \times \pi_{t-h}(u)du dt  \}\nonumber\\
		=&  \min_{\pi \in \mathcal{P}( \mathbb{R}^m ) } \frac{1}{2}  E \{  \int_{0}^{\infty } e^{-\rho t} \int_{\mathbb{R}^m}  [	
		(u -K\hat{X}^{\pi}_{t\mid t-h})'\Omega \nonumber \\
		&\times  (u -K\hat{X}^{\pi}_{t\mid t-h}) ] \pi_{t-h}(u)du  dt \}+   E \{  \int_{0}^{\infty } e^{-\rho t}  \nonumber\\
		& \times  \int_{\mathbb{R}^m} \lambda  \ln \pi_{t-h}(u) \pi_{t-h}(u)du dt  \}  
		+ \frac{1}{2}  E \{   x'Px \nonumber \\
		&  -x' \int_{0}^{h}\Lambda(s)\hat{x}_{0\mid s-h}ds \} .
	\end{align}
	Note that $\pi \in \mathcal{P}( \mathbb{R}^m )$, if and only if
	\begin{align*}
		\int_{ \mathbb{R}^m }\pi_{t-h}(u)du = 1 \ and \  \pi_{t-h}(u)\ge0 \ a.e.\ on \  { \mathbb{R}^m }.
	\end{align*}	
	%
	%
	Solving the constrained minimization problem \eqref{reypi}, we obtain that 
	\begin{align}\label{entropy pi}
		\pi_{t-h}^*  
		= &\frac{1}{N_\pi}  exp\{   -\frac{1}{2}  [u-K\hat{X}^{\pi}_{t\mid t-h} ]'  \frac{1}{\lambda}(R+\bar{B}' P\bar{B} ) \nonumber
		\\
		& \times [u-K\hat{X}^{\pi}_{t\mid t-h}  ] \} ,
	\end{align}
	with the normalization constant
	\begin{align*}
		N_\pi = & \int_{U}     exp\{   -\frac{1}{2}   [u-K\hat{X}^{\pi}_{t\mid t-h} ]'  \frac{1}{\lambda}(R+\bar{B}' P\bar{B} ) \nonumber
		\\
		& \times [u-K\hat{X}^{\pi}_{t\mid t-h}] \}  \pi_{t-h}   du,
	\end{align*} 
	which indicates that the optimal control follows a Gaussian distribution, where the mean and variance of the optimal policy \(\pi_{t-h}^*\) are presented in \eqref{mean}–\eqref{var}.
	
	Submitting \eqref{entropy pi} into \eqref{reypi} yields
	\begin{align*}		
		V&^{\pi}(x)\nonumber \\
		=&\min_{\pi \in \mathcal{P}( \mathbb{R}^m ) } \frac{1}{2}  E \{  \int_{0}^{\infty }  e^{-\rho t } \int_{\mathbb{R}^m}  [	
		(u-K\hat{X}^{\pi}_{t \mid t -h})'\Omega(u -K \nonumber \\
		& \times  \hat{X}^{\pi}_{t \mid t-h}) ] \mathcal{N}(u \mid  \mu^{*}_{t-h} , \Sigma^{*}_{t-h})  du ds\}   -  \frac{\lambda }{2\rho} \ln((2\pi e)^m   \nonumber\\
		&\times  \left |\lambda\Omega^{-1}   \right | )      + \frac{1}{2}  E \{   [ x'Px -x' \int_{0}^{h}\Lambda(s) \hat{x}_{0\mid s-h}ds] \} \nonumber \\	
		=	& \frac{1}{2}  \{  x'Px-x'\int_{0}^{h}\Lambda(s)\hat{x}^{\pi}_{0\mid s-h}ds\}	-  \frac{\lambda }{2\rho} [ \ln((2\pi e)^m  \nonumber \\
		& \times  \left |  \lambda\Omega^{-1}   \right | )  -1    ] .
	\end{align*}
	
	By Theorem 6 in \cite{Wang2019} and the admissibility of \(u^*(t-h)\) presented in part \textbf{(i)}, we conclude that \(\pi^*_{t-h}\) given by \eqref{pi*} is admissible.
	Accordingly, the admissible optimal control and value function of Problem \ref{problem 2} are shown in \eqref{pi*} and \eqref{Vpifunction}.
	The proof is now completed.
	\hfill $\blacksquare$


\begin{thebibliography}{99}
		\bibitem{finance}
		L. Chen, Z. Wu, Maximum principle for the stochastic optimal control problem with delay and application, Automatica, Vol. 46, no. 6, pp. 1074-1080, June 2010.
		
		\bibitem{network}
		J. Baillieul and P. J. Antsaklis, Control and Communication Challenges in Networked Real-Time Systems,  Proceedings of the IEEE, vol. 95, no. 1, pp. 9-28, Jan. 2007.
		
		\bibitem{power system}
		Z. Wang, J. Dai, H. Zhang and J. Zhang, Observer-Based Finite-Time Fuzzy Load Frequency Control for Multiarea Nonlinear Power Systems Under Input Delays and Cyber Attacks, IEEE Transactions on Smart Grid, vol. 16, no. 4, pp. 3336-3345, July 2025.
		
		\bibitem{stable1}
		S. I. Niculescu, Delay effects on stability: a robust control approach, Springer London, 2002.
		
		\bibitem{stalbe2}
		N. Bekiaris-Liberis and M. Krstic, Nonlinear control under nonconstant delays, Society for Industrial and Applied Mathematics, 2013.
		
		\bibitem{delay1}
		Z. Artstein, Linear systems with delayed controls: A reduction, IEEE Transactions on Automatic Control, vol. 27, no. 4, pp. 869-879, August 1982.
		
		\bibitem{delay4}
		C. Hua, Q. Wang and X. Guan, Adaptive Tracking Controller Design of Nonlinear Systems With Time Delays and Unknown Dead-Zone Input, IEEE Transactions on Automatic Control, vol. 53, no. 7, pp. 1753-1759, Aug. 2008.
		
		\bibitem{delay2}
		N. Espitia and W. Perruquetti, Predictor-Feedback Prescribed-Time Stabilization of LTI Systems With Input Delay, IEEE Transactions on Automatic Control, vol. 67, no. 6, pp. 2784-2799, June 2022.
		
		\bibitem{delay3}
		A. Wu, S. Shen, J. Zhang and J. Mei, A Model Reduction Approach for Discrete-Time Coupled Systems With Input Delays via Bivariate Fundamental Matrices, IEEE Transactions on Automatic Control, vol. 71, no. 5, pp. 2950-2965, May 2026.
		
		\bibitem{xie}
		X. Chen and X. Xie, SARF Copilot: A Holistic Safety-Aware Resilient Fuzzy Control Synthesis for Human–Machine Shared Steering System, IEEE Transactions on Industrial Electronics, vol. 73, no. 6, pp. 9354-9366, June 2026.	
		
		\bibitem{wangwenjing}
		W. Wang, J. Xu, H. Zhang, M. Fu, Exact controllability of rational expectations model with multiplicative noise and input delay, Journal of Automation and Intelligence, Vol. 3, no. 1, pp. 19-25, January 2024.	
		
		
		\bibitem{stoch1}
		B. Larssen, Dynamic programming in stochastic control of systems with delay. Stochastics and Stochastic Reports, vol. 74, pp. 651–673, 2002.
		
		\bibitem{stoch2}
		Y. Shen, Q. Meng, P. Shi, Maximum principle for mean-field jump–diffusion stochastic delay differential equations and its application to finance, Automatica, vol. 50, no. 6, pp. 1565-1579, June 2014.
		
		\bibitem{stoch3}
		H. Zhang and J. Xu, Control for Itô Stochastic Systems With Input Delay, IEEE Transactions on Automatic Control, vol. 62, no. 1, pp. 350-365, Jan. 2017.
		
		\bibitem{stoch4}
		X. Xue, J. Xu and H. Zhang, Stabilization of First-Order Hyperbolic PDE Systems With Input Delay and Multiplicative Noise, IEEE Transactions on Automatic Control, vol. 71, no. 5, pp. 3372-3379, May 2026.
		
		
		
		\bibitem{Wang2019}
		H. Wang, T. Zariphopoulou and  X. Zhou, Exploration versus Exploitation in Reinforcement Learning: A Stochastic Control Approach, 2019, arXiv preprint arXiv:1812.01552.
		
		\bibitem{Wang2020}
		H. Wang and  X. Zhou, Continuous-time mean–variance portfolio selection: A reinforcement learning framework, Mathematical Finance, vol. 30, no. 4, pp. 1273-1308, June 2020.
		
		\bibitem{rl1}
		Y. Jia and X. Y. Zhou, Q-learning in continuous time, Journal of Machine Learning Research, vol. 24, no. 1, pp. 1-61, Jan. 2023.
		
		\bibitem{rl2}
		L. Szpruch, T. Treetanthiploet and Y. Zhang, Optimal Scheduling of Entropy Regularizer for Continuous-Time Linear-Quadratic Reinforcement Learning, SIAM Journal on Control and Optimization, vol. 62, no. 1, pp. 135-166, Jan. 2024.
		
		
		\bibitem{rl3} Z. Chen, Q. Zhang, Backward Stochastic Control System with Entropy Regularization, SIAM Journal on Control and Optimization, vol. 63, no. 3, pp. 1981-2006, June 2025.
		
		\bibitem{11}
		J. Xu and H. Zhang, Exponential Mean Square Stabilization for Itô Stochastic Systems with Input Delay, 2019 Chinese Control Conference (CCC), Guangzhou, China, 2019, pp. 1405-1410.
		
		
		\bibitem{HJB1}
		F. De Feo, S. Federico and A. Świech, Optimal Control of Stochastic Delay Differential Equations and Applications to Path-Dependent Financial and Economic Models, SIAM Journal on Control and Optimization, vol. 62, no. 3, pp. 1490-1520, May 2024.
		
		\bibitem{HJB2}
		W. Fleming and H. Soner, Controlled Markov Processes and Viscosity Solutions, New York, NY, USA: Springer, 2006.
		
		\bibitem{Lyapunov stability theorem}
		R. Kalman and J. Bertram, Control system analysis and design via the second method of lyapunov: (I) continuous-time systems (II) discrete time systems, IRE Transactions on Automatic Control, vol. 4, no. 3, pp. 112-112, December 1959.
		
		
		
		
		
		
	\end{thebibliography}
\end{document}